\documentclass[a4paper,11pt]{amsart}
\usepackage[T1]{fontenc}
\usepackage[english]{babel}
\usepackage[cp1252]{inputenc}
\usepackage{amsthm}
\usepackage{amsmath}
\usepackage{amsfonts}
\usepackage{stmaryrd}
\usepackage{amssymb}
\usepackage{enumitem}
\newcommand{\Ca}{\mathrm{card\,}}
\newcommand{\e}{\varepsilon}
\newcommand{\di}{\mathrm{d}}

\newcommand{\alp}{\omega(P)+2}

\newtheorem{theorem}{Theorem}[section]

\newtheorem{proposition}{Proposition}[section]

\begin{document}

\address{NWF I-Mathematik, Universit\"at Regensburg, 93040 Regensburg, Germany.}
\email{\href{mailto:farid.madani@mathematik.uni-regensburg.de}{Farid.Madani@mathematik.uni-regensburg.de}}

\subjclass{}

\keywords{}
\author{Farid Madani}
\title{Hebey--Vaugon conjecture II}
\date{}
\maketitle

\section{Introduction} 
Let $(M,g)$ be a compact Riemannian manifold of dimension $n\geq 3$. Denote by $I(M,g)$, $C(M,g)$ and $R_g$ the isometry group, the conformal transformations group and the scalar curvature, respectively. Let $G$ be a subgroup of the isometry group $I(M,g)$. The equivariant Yamabe problem can be formulated as follows: in the conformal class of $g$, there exists a $G-$invariant metric with constant scalar curvature. Assuming the positive mass theorem and the Weyl vanishing conjecture  (for more details on the subject, see \cite{KMS}, \cite{Mar} and the references therein), E.~Hebey and M.~Vaugon~\cite{HV} proved that this problem has  solutions. Moreover, they proved that the infimum of Yamabe functional 
\begin{equation}\label{Yamabe func}
I_g(\varphi)=\frac{\int_M |\nabla\varphi|^2+\frac{n-2}{4(n-1)}R_g\varphi^2\di v}{\|\varphi\|^2_{\frac{2n}{n-2}}}
\end{equation}
over $G-$invariant nonnegative functions is achieved by a smooth positive $G-$invariant function. This function is a solution of the Yamabe equation, which is the Euler--Lagrange equation of $I_g$:
\begin{equation*}
\frac{4(n-1)}{n-2}\Delta_g \varphi+ R_g\varphi= \mu\varphi^{\frac{n+2}{n-2}}
\end{equation*}

One of the consequences of these results is  that the following conjecture due to Lichnerowicz \cite{Lic} is true.
\vspace{0.3cm}
\paragraph{\textsc{Lichnerowicz conjecture}}\emph{ For every compact Riemannian manifold $(M,g)$ which is not conformal to the unit sphere $S^n$ endowed with its standard metric $g_s$, there exists a metric $\tilde g$ conformal to $g$ for which $I(M,\tilde g)=C(M,g)$, and the scalar curvature $R_{\tilde g}$ is constant.}\\

The classical Yamabe problem, which consists of finding a conformal metric with constant scalar curvature on a compact Riemannian manifold, is a particular case of the equivariant Yamabe problem (it corresponds to $G=\{\mathrm{id}\}$). This problem was completely solved by H.~Yamabe~\cite{Yam}, N.~Trudinger~\cite{Trud}, T.~Aubin~\cite{Aub} and R.~Schoen~\cite{Schoen}.  The main idea to prove the existence of positive minimizers for $I_g$ is  to show that if $(M,g)$ is not conformal to the sphere endowed with its standard metric, then
\begin{equation}\label{sss}
\mu(g):=\inf_{C^\infty(M)} I_g(\varphi)<\frac{1}{4}n(n-2)\omega_n^{2/n}
\end{equation}
where $\omega_n$ is the volume of the unit sphere $S^n$.

T.~Aubin~\cite{Aub} proved \eqref{sss} in some cases by constructing a test function $u_\e$ satisfying:
\begin{equation*}
I_g(u_\e)<\frac{1}{4}n(n-2)\omega_n^{2/n}
\end{equation*}

He conjectured that \eqref{sss} always holds except for the sphere. R.~Schoen  constructed another test function which involves the Green function of the conformal Laplacian $\Delta_g+\frac{n-2}{4(n-1)}R_g$. Using the positive mass theorem, R.~Schoen proved~\eqref{sss} for all compact manifolds which are not conformal to $(S^n,g_s)$. The solution of  the Yamabe problem follows. \\

Later, E.~Hebey and M.~Vaugon~\cite{HV} showed that we can generalize \eqref{sss} for the equivariant case as follows:\\
Denote by $O_G(P)$ the orbit of $P\in M$ under $G$ and by $\Ca O_G(P)$ its cardinal. Let $C^\infty_G(M)$ be the set of smooth $G$-invariant functions and 
\begin{equation*}
\mu_G(g):=\inf_{ C_G^\infty(M)} I_g(\varphi) 
\end{equation*}

Following E.~Hebey and M.~Vaugon~\cite{HV4,HV}, we define the integer $\omega(P)$ at a point $P$ as 
$$\omega(P)=\inf \{i\in \mathbb N/\|\nabla^i W_g(P)\|\neq 0\}\; (\omega(P)=+\infty\text{ if }\forall i\in \mathbb N,\;\|\nabla^i W_g(P)\|=0)$$

\medskip

\paragraph{\textsc{Hebey--Vaugon conjecture}}\emph{Let $(M,g)$ be a compact Riemannian manifold of dimension $n\geq 3$ and $G$ be a subgroup of $I(M,g)$. If $(M,g)$ is not conformal to $(S^n, g_{s})$ or if the action of $G$ has no fixed point, then the following inequality holds }
\begin{equation}\label{HVI}
\mu_G(g)<\frac{1}{4}n(n-2)\omega_n^{2/n}(\inf_{Q\in M}\Ca O_G(Q))^{2/n}
\end{equation}

E.~Hebey and M.~Vaugon showed that if this conjecture holds, then it implies that the equivariant Yamabe problem has minimizing solutions and the Lichnerowicz conjecture is also true. Notice that if $G=\{\mathrm{id}\}$, then this conjecture corresponds to \eqref{sss}.\\
Let us recall the results already known about this conjecture. Assuming the positive mass theorem, E.~Hebey and M.~Vaugon \cite{HV} proved the following:
\begin{theorem}[E.~Hebey and M.~Vaugon]\label{HV theorem}
Let $(M,g)$ be a smooth compact Riemannian manifold of dimension $n\geq 3$  and $G$ be a subgroup of $I(M,g)$.
We always have :  
$$\mu_G(g)\leq \frac{1}{4}n(n-2)\omega_n^{2/n}(\inf_{Q\in M}\Ca O_G(Q))^{2/n}$$
 and inequality \eqref{HVI} holds if at least one of the following conditions is satisfied. 
 \begin{enumerate}[label=\arabic*.]
 \item The action of $G$  on $M$ is free.
 \item $3\leq \dim M\leq 11$.
 \item\label{item 3}There exists a point $P\in M$ with finite minimal orbit  under $G$ such that $\omega(P)>(n-6)/2$ or $\omega(P)\in \{0,1,2\}$.
 \end{enumerate}
\end{theorem}
We have also the following result obtained by the author in \cite{Mad2}: 
\begin{theorem}\label{cor prin}
Hebey--Vaugon conjecture  holds  for every smooth compact Riemannian manifold $(M,g)$ of dimension $n\leq 37$.
\end{theorem} 

The main result of this paper is  the following:

\begin{theorem}\label{main}
If there exists a point $P\in M$ such that $\omega(P)\leq(n-6)/2$, then 
\begin{equation}\label{for p}
\mu_G(g)<\frac{1}{4}n(n-2)\omega_n^{2/n}(\Ca O_G(P))^{2/n}
\end{equation}
\end{theorem}

Note that if we assume the positive mass theorem, then Theorem \ref{main} and Theorem \ref{HV theorem} implies that Hebey--Vaugon conjecture holds.\\ 
The proof of Theorem \ref{main} doesn't require the positive mass theorem. If $\Ca O_G(Q)=+\infty$ for all $Q\in M$, then \eqref{HVI} holds. So we have to consider only the case when there exists a point in $M$ with finite orbit.
From now until the end of this paper, we suppose that $P\in M$ is contained in a finite orbit and  $\omega(P)\leq \frac{n-6}{2}$. The assumption $\omega(P)\leq \frac{n-6}{2}$ deletes the case $(M,g)$ is conformal to  $(S^n,g_s)$. 

\section{G-invariant test function}

In order to prove Theorem \ref{main} and \ref{cor prin}, we construct from the function $\varphi_{\e,P}$ defined below a $G$-invariant test function $\phi_\e$ such that
\begin{equation}\label{a prouver}
I_g(\phi_\e)<\frac{1}{4}n(n-2)\omega_{n}^{2/n}(\Ca O_G(P))^{2/n}
\end{equation}

Let us recall the construction in \cite{Mad2} of $\varphi_{\e,P}$. Let $\{x^j\}$ be the geodesic normal coordinates  in the neighborhood of $P$ and define $r=|x|$ and $\xi^j=x^j/r$. Without loss of generality, we suppose that $\det g=1+O(r^N)$, with $N>0$ sufficiently large (for the existence of  such coordinates for a $G-$invariant conformal class, see \cite{HV}, \cite{LP}).

\begin{equation*}\label{phiepsilon}
\varphi_{\e,P}(Q)=(1- r^{\alp} f(\xi))u_{\e,P}(Q)
\end{equation*}

\begin{equation*}\label{uepsilon}
u_{\e,P}(Q)=\begin{cases}\biggl(\displaystyle\frac{\varepsilon}{r^2+\varepsilon^2}\biggr)^{\frac{n-2}{2}}-\biggl(\frac{\varepsilon}{\delta^2+\varepsilon^2}\biggr)^{\frac{n-2}{2}} &\mbox{ if }Q\in B_{P}(\delta)\\
0 &\mbox{ if }Q\in M-B_{P}(\delta)
\end{cases}
\end{equation*}
for all $Q\in M$, where $r=d(Q,P)$ is the distance between $P$ and $Q$, and $B_{P}(\delta)$ is the geodesic ball of center $P$ and radius $\delta$ fixed sufficiently small. $f$ is a function  depending only on $\xi$ (defined on $S^{n-1}$), chosen such that $\int_{S^{n-1}}f d\sigma=0$. \\ 

Let $\bar R$ be the leading part in the Taylor expansion of the scalar curvature $R_g$ in a neighborhood of $P$ and $\mu(P)$ is its degree. Hence,
$$R_g(Q)=\bar R +O(r^{\mu(P)+1})$$  
$$\bar R=r^{\mu(P)}\sum_{|\beta|= \mu(P)}\nabla_\beta R_g(P)\xi^\beta$$ 
We summarize some properties of $\bar R$ in the following proposition.
\begin{proposition}\label{propo R bar}
 \begin{enumerate}[label=\arabic*.]
  \item $\bar R$ is a homogeneous polynomial of degree $\mu(P)$ and is invariant under the action of the stabilizer group of $P$.
  \item We always have  $\mu(P)\geq \omega(P)$ 
  \item if $\mu(P)\geq \omega(P)+1$, then $\int_{S^{n-1}(r)}R\di\sigma<0$ for $r>0$ sufficiently small.
  \item \label{item c}If $\mu(P)=\omega(P)$, then there exist eigenfunctions $\varphi_k$ of the Laplacian on $S^{n-1}$ such that the restriction of $\bar R$ to the sphere is given by 
$$\bar R|_{S^{n-1}}=\sum_{k=1}^{q}\nu_k\varphi_{k}$$
where $q\leq [\omega(P)/2]$, $\Delta_s\varphi_k=\nu_k\varphi_k$ and $\nu_k=(\omega-2k+2)(n+\omega-2k)$ are the eigenvalues of $\Delta_s$ with respect to the standard metric $g_s$ of $S^{n-1}$.
 \end{enumerate}

\end{proposition}

Since the scalar curvature is invariant under the action of the isometry group $I(M,g)$, $\bar R$ is invariant under the action of the stabilizer of $P$. The second statement of Proposition \ref{propo R bar} is proven by E.~Hebey and M.~Vaugon (\cite{HV}, section 8)  and the third one by T.~Aubin~(\cite{Aub3}, section 3). In this case the conjecture holds immediately, by choosing  $f=0$, $\varphi_{\e,P}=u_{\e,P}$ (see \cite{Mad3,Mad2} for more details). \\
From now we suppose that $\mu(P)=\omega(P)$. Using the fact that $\bar R$ is homogeneous polynomial of degree $\omega (P)$ and the fact that for all $j\leq \omega(P)-1$
\begin{equation}\label{identity}
|\nabla^jR_g(P)|=0,\;\;\Delta_g^{j+2}R_g(P)=0\text{ and }|\nabla\Delta_g^{j+2}R_g(P)|=0
\end{equation}
we deduce that $\Delta_{\mathcal E}^{[\omega(P)/2]}\bar R= 0$. Hence, if we restrict $\bar R$ to the sphere, we get the decomposition of item \ref{item c} in Proposition \ref{propo R bar}. The proof of \eqref{identity} is given in \cite{HV}, section 8.

Using the split of $\bar R$  given in Proposition \ref{propo R bar}, we proved in \cite{Mad2} that  if the cardinal of $O_G(P)$ is minimal and $\omega(P)\leq 15$ , then there exists $c\in\mathbb R$ such that for $f=c\bar R\vert_{S^{n-1}}$, the function 
$$\phi_\e=\sum_{P_i\in O_G(P)}\varphi_{\e,P_i}$$
is $G$-invariant and satisfies \eqref{a prouver}, which proves Theorem \ref{cor prin}. Moreover, we proved the following theorem:

\begin{theorem}\label{ddd}
If $\omega(P)\leq (n-6)/2$, then there exist $c_k\in\mathbb R$, such that for $f=\sum_{k=1}^qc_k\varphi_k$, the function $\varphi_{\e, P}$ satisfies
\begin{equation}\label{inega omega}
I_g(\varphi_{\e,P})<\frac{1}{4}n(n-2)\omega_{n}^{2/n}
\end{equation}
\end{theorem}
The proof of Theorem \ref{ddd} is technical and uses Proposition \ref{propo R bar}. It is given in \cite{Mad2} (see also \cite{Mad3} for a detailed proof).
  
Below, we show that using Theorem \ref{ddd}, we can construct a $G-$invariant function $\phi_\e$ which satisfies \eqref{a prouver} for $\omega(P)\leq \frac{n-6}{2}$ (the cardinal of $O_G(P)$ is not necessarily minimal). It implies Theorem \ref{main}.
\begin{proof}[Proof of Theorem \ref{main}]
Let $H\subset G$ be the stabilizer of $P$. We consider the function $f=\sum_{k=1}^qc_k\varphi_k$ of Theorem \ref{ddd}.  Using the exponential map on $P$ as a local chart, we can view $f$ and $\varphi_k$ as functions defined over the unit sphere of $T_PM$, the tangent space of $M$ on $P$. Let $h$ be an isometry in $H$. 
$$h_*(P):(T_PM,g_P)\rightarrow (T_PM,g_P)$$ 
is the linear tangent map of $h$ on $P$. It is a linear isometry with respect to the inner product $g_P$ which is Euclidean. $h_*(P)$ conserves the unit sphere $S^{n-1}\subset T_PM$ and the Laplacian. We already know that the function $\bar R=r^{\omega(P)}\sum_{k=1}^{q}\nu_k\varphi_{k}$ is $H$-invariant. Notice that $\varphi_k$ and $\varphi_j$ belong to two different eigenspaces if $k\neq j$. Since, isometries conserve the Laplacian and $\varphi_k$ are eigenfunctions  of the Laplacian on the sphere endowed with its standard metric, it yields that $\varphi_k$ and $f$ are $H$-invariant. On the other hand, we have the following bijective map: 
\begin{equation*}
 \begin{split}
 G/H &\longrightarrow O_G(P)\\
\sigma H & \longmapsto \sigma(P)
 \end{split}
\end{equation*} 
Since $f$ is $H$-invariant, $\varphi_{\e,P}$ is $H-$invariant and the function 
$$\phi_\e=\sum_{\sigma\in G/H}\varphi_{\e,P}\circ\sigma^{-1}$$ 
is $G-$invariant and satisfies \eqref{a prouver}. 
\end{proof}
\section*{Acknowledgments}
The author would like to thank Bernd Ammann for his remarks  and the helpful discussions.  He also thanks Emmanuel Hebey for his suggestions.

\bibliographystyle{amsplain}
\bibliography{bibliographie}

\end{document}